\documentclass[9pt]{article}
\usepackage{mathrsfs}
\usepackage{amsthm}
\usepackage{amssymb}
\usepackage{amsmath}
\usepackage{graphicx}
\usepackage{color}
\usepackage{amsfonts}
\usepackage{float}
\usepackage{cite}
\usepackage [latin1]{inputenc}
\usepackage[text={140mm,210mm},left=35mm,vmarginratio=1:1]{geometry}
\newtheorem{theorem}{Theorem}[section]
\newtheorem{remark}{Remark}[section]
\newtheorem{lemma}[theorem]{Lemma}

\newtheorem{corollary}[theorem]{Corollary}

\numberwithin{equation}{section}
\normalsize

\begin{document}
\title{\textbf{Equilibrium moderate deviations for occupation times of SSEP on regula trees}}

\author{Xiaofeng Xue \thanks{\textbf{E-mail}: xfxue@bjtu.edu.cn \textbf{Address}: School of Mathematics and Statistics, Beijing Jiaotong University, Beijing 100044, China.}\\ Beijing Jiaotong University}

\date{}
\maketitle

\noindent {\bf Abstract:} In this paper, we are concerned with the symmetric simple exclusion process on the regula tree $\mathbb{T}^d$ for $d\geq 2$. Our main result gives moderate deviation principles of occupation times of the process starting from an invariant product measure. Two replacement lemmas play key roles in the proof of our main result. To obtain these replacement lemmas, we utilize duality relationships between the symmetric exclusion process and two types of random walks on $\mathbb{T}^d$ and $\left(\mathbb{T}^d\right)^2$ respectively.

\quad


\section{Introduction}\label{section one}
\subsection{The model}\label{subsection 1.1}
In this paper, we are concerned with the symmetric simple exclusion process (SSEP) on a regular tree $\mathbb{T}^d$ with $d\geq 2$, where each vertex has $d+1$ neighbors. For later use, for any $x,y\in \mathbb{T}^d$, we write $x\sim y$ when they are neighbors with each other. The SSEP $\{\eta_t\}_{t\geq 0}$ is a continuous-time Markov process with state space $\{0,1\}^{\mathbb{T}^d}$, i.e., each vertex on $\mathbb{T}^d$ is occupied by a particle or vacant. The generator $\mathcal{L}$ of $\{\eta_t\}_{t\geq 0}$ is given by
\begin{equation}\label{equ 1.1. generator}
\mathcal{L}f(\eta)=\frac{1}{2}\sum_{x\in \mathbb{T}^d}\sum_{y:y\sim x}\left[f(\eta^{x,y})-f(\eta)\right]
\end{equation}
for any $\eta\in \{0, 1\}^{\mathbb{T}^d}$ and $f$ from $\mathbb{T}^d$ to $\mathbb{R}$ depending on finite coordinates, where
\[
\eta^{x,y}(z)=
\begin{cases}
\eta(z) & \text{~if~}z\neq x\text{~and~}z\neq y,\\
\eta(y) & \text{~if~}z=x,\\
\eta(x) & \text{~if~}z=y.
\end{cases}
\]
According to the definition of $\mathcal{L}$, in the SSEP, all particles perform simple random walks on $\mathbb{T}^d$, where a particle jumps from a vertex $x$ to each neighbor $y$ of $x$ at rate $1$. However, any jump to an occupied vertex is suppressed, since on each vertex there is at most one particle. For a detailed survey of the exclusion process, see Chapter 8 of \cite{Lig1985} and Part ${\rm \uppercase\expandafter{\romannumeral3}}$ of \cite{Lig1999}.

For given $0<p<1$, we denote by $\nu_p$ the product measure on $\mathbb{T}^d$ under which $\{\eta(x)\}_{x\in \mathbb{T}^d}$ are independent and
\[
\nu_p\left(\eta(x)=1\right)=p=1-\nu_p\left(\eta(x)=0\right)
\]
for all $x\in \mathbb{T}^d$. We denote by $\mathbb{P}_{\nu_p}$ the probability measure of the SSEP $\{\eta_t\}_{t\geq 0}$ starting from $\nu_p$. We denote by $\mathbb{E}_{\nu_p}$ the expectation with respect to $\mathbb{P}_{\nu_p}$. According to the definition of $\mathcal{L}$, it is easy to check that
\begin{equation}\label{equ 1.2 equilibrium}
\int f(\eta)\mathcal{L}g(\eta)\nu_p(d\eta)=\int g(\eta)\mathcal{L}f(\eta)\nu_p(\eta)
\end{equation}
for any $f,g$ from $\{0,1\}^{\mathbb{T}^d}$ to $\mathbb{R}$ depending on finite coordinates. As a result, $\nu_p$ is a reversible measure of the SSEP $\{\eta_t\}_{t\geq 0}$. In this paper, we investigate the moderate deviation principle of the occupation time of the SSEP starting from $\nu_p$. For mathematical details, see Section \ref{section two}.

\subsection{Duality relationships}\label{subsection 1.2}
For later use, in this section we recall the duality relationship between SSEP and random walks. We denote by $\{V_t\}_{t\geq 0}$ the continuous-time simple random walk on $\mathbb{T}^d$ with generator $\Omega_1$ given by
\[
\Omega_1h(x)=\sum_{y:y\sim x}\left(h(y)-h(x)\right)
\]
for any $x\in \mathbb{T}^d$ and bounded $h$ from $\mathbb{T}^d$ to $\mathbb{R}$. For any $t\geq 0$ and $x,z\in \mathbb{T}^d$, we denote by $p_t(x,z)$ the probability $\mathbb{P}\big(V_t=z\big|V_0=x\big)$. We write $V_t$ as $V_t^x$ when $V_0=x$. We denote by $\{Y_t\}_{t\geq 0}$ the continuous-time random walk on
\[
\mathbb{Y}:=\left(\mathbb{T}^d\right)^2\setminus \left\{(y,w)\in \left(\mathbb{T}^d\right)^2:y=w\right\}
\]
with generator $\Omega_2$ given by
\[
\Omega_2g(x,z)=\sum_{(y,w)\in \mathbb{Y}}Q\left((x,z), (y,w)\right)\left(g(y,w)-g(x,z)\right)
\]
for any $(x,z)\in \mathbb{Y}$ and bounded $g$ from $\mathbb{T}^d$ to $\mathbb{R}$, where
\[
Q\left((x,z), (y,w)\right)=
\begin{cases}
1 & \text{~if~}x\not\sim z, y\sim x \text{~and~} w=z,\\
1 & \text{~if~}x\not\sim z, y=x\text{~and~} w\sim z,\\
1 & \text{~if~}x\sim z, y=z \text{~and~}w=x,\\
1 & \text{~if~}x\sim z, y\sim x, y\neq z\text{~and~}w=z,\\
1 & \text{~if~}x\sim z, w\sim z, w\neq x\text{~and~}y=x,\\
0 & \text{~else}.
\end{cases}
\]
For $k=1,2$, we denote by $Y_t(k)$ the $k$th component of $Y_t$. According to the expression of $Q(\cdot, \cdot)$, $\{Y_t(1)\}_{t\geq 0}$ and $\{Y_t(2)\}_{t\geq 0}$ perform independent simple random walks on $\mathbb{T}^d$ except that they exchange positions with each other at rate $1$ when they are neighbors to avoid collision. For any $(x,z), (y,w)\in \mathbb{Y}$, we denote by $q_t\left((x,z), (y,w)\right)$ the probability $\mathbb{P}\left(Y_t=(y,w)\big|Y_0=(x,z)\right)$.

For any $\eta\in \{0,1\}^{\mathbb{T}^d}$, we denote by $\mathbb{P}_{\eta}$ the probability measure of the SSEP starting from $\eta$ and by $\mathbb{E}_{\eta}$ the expectation with respect to $\mathbb{P}_{\eta}$. We have the following duality relationship.

\begin{lemma}\label{lemma dual}
For any $t\geq 0$, $h_1$ from $\{0, 1\}$ to $\mathbb{R}$, $h_2$ from $\{0,1\}^2$ to $\mathbb{R}$ and $x, y\in \mathbb{T}^d$ such that $x\neq y$,
\begin{equation}\label{equ 1.3 1vertexdual}
\mathbb{E}_\eta h_1(\eta_t(x))=\mathbb{E}h_1(\eta(V_t^x))=\sum_{z\in \mathbb{T}^d}p_t(x,z)h_1(\eta(z))
\end{equation}
and
\begin{align}\label{equ 1.4 2vertexdual}
\mathbb{E}_\eta h_2(\eta_t(x), \eta_t(y))&=\mathbb{E}h_2\left(\eta\left(Y_t^{(x,y)}(1)\right),\eta\left(Y_t^{(x,y)}(2)\right)\right)\notag \\
&=\sum_{(z,w)\in \mathbb{Y}}q_t\left((x,y), (z,w)\right)h_2\left(\eta(z), \eta(w)\right).
\end{align}

\end{lemma}

\proof[Proof of Lemma \ref{lemma dual}]

We only prove Equation \eqref{equ 1.4 2vertexdual} since Equation \eqref{equ 1.3 1vertexdual} follows from a similar analysis. For any $(x,y)\in \mathbb{Y}$, we supplementarily define
\[
Q\left((x,y), (x,y)\right)=-\sum_{(z,w):(z,w)\neq (x,y)}Q\left((x,y), (z,w)\right).
\]
For any $t\geq 0$ and $x,y\in \mathbb{Y}$, we define
\[
H_t(x,y)=\mathbb{E}_\eta h_2\left(\eta_t(x), \eta_t(y)\right).
\]
According to the definition of $\mathcal{L}$ and Kolmogorov-Chapman Equation, if $x\not\sim y$, then
\begin{align*}
\frac{d}{dt}H_t(x,y)=&\sum_{z:z\sim x}\Big(\mathbb{E}_\eta h_2\left(\eta_t(z), \eta_t(y)\right)-\mathbb{E}_\eta h_2\left(\eta_t(x), \eta_t(y)\right)\Big)\\
&+\sum_{w:w\sim y}\Big(\mathbb{E}_\eta h_2\left(\eta_t(x), \eta_t(w)\right)-\mathbb{E}_\eta h_2\left(\eta_t(x), \eta_t(y)\right)\Big).
\end{align*}
If $x\sim y$, then
\begin{align*}
&\frac{d}{dt}H_t(x,y)\\
&=\mathbb{E}_{\eta}h_2\left(\eta_t(y), \eta_t(x)\right)-\mathbb{E}_{\eta}h_2\left(\eta_t(x), \eta_t(y)\right)\\
&\text{\quad}+\sum_{z:z\sim x,z\neq y}\Big(\mathbb{E}_\eta h_2\left(\eta_t(z), \eta_t(y)\right)-\mathbb{E}_\eta h_2\left(\eta_t(x), \eta_t(y)\right)\Big)\\
&\text{\quad}+\sum_{w:w\sim y, w\neq x}\Big(\mathbb{E}_\eta h_2\left(\eta_t(x), \eta_t(w)\right)-\mathbb{E}_\eta h_2\left(\eta_t(x), \eta_t(y)\right)\Big).
\end{align*}
In conclusion, $\frac{d}{dt}H_t=QH_t$ and hence Equation \eqref{equ 1.4 2vertexdual} holds according to the initial condition where $H_0(z,w)=h_2\left(\eta(z), \eta(w)\right)$ for all $(z,w)\in \mathbb{Y}$.

\qed

\subsection{Occupation times}\label{subsection 1.3}
In this subsection we recall the definition of occupation times. For any $x\in \mathbb{T}^d$ and $t\geq 0$, the occupation time $X_t^x$ of $\{\eta_s\}_{0\leq s\leq t}$ on $x$ is defined as
\[
X_t^x=\int_0^t\eta_s(x)ds.
\]
For any $0<p<1$, by Equations \eqref{equ 1.2 equilibrium}, we have
\[
\mathbb{P}_{\nu_p}\left(\eta_t(x)=1\right)=p
\]
for all $x\in \mathbb{T}^d$ and $t\geq 0$. Hence, when $\{\eta_t\}_{t\geq 0}$ starts from $\nu_p$, the centered occupation time $\xi_t^x$ on $x$ is naturally defined as
\[
\xi_t^x=\int_0^t\left(\eta_s(x)-p\right)ds.
\]
Since 1980s, limit theorems of occupation times of SSEP are popular research topics. Reference \cite{Kipnis1987} investigates central limit theorems of occupation times of SSEP on lattices $\mathbb{Z}^d$. It is shown in \cite{Kipnis1987} that, when $d\geq 1$ and the SSEP on $\mathbb{Z}^d$ starts from $\nu_p$ (lattice version), there exists $n_t=n_t(d)$ such that $\frac{1}{n_t}\xi_t^x$ (lattice version) converges weakly to some Gaussian random variable $\xi$ as $t\rightarrow+\infty$. In detail, $n_t$ performs following dimension-dependent phase transition. For $d\geq 1$,
\[
n_t(d)=
\begin{cases}
t^{3/4} & \text{~if~}d=1,\\
\sqrt{t\log t} & \text{~if~}d=2,\\
t^{1/2} & \text{~if~}d\geq 3.
\end{cases}
\]
Large and moderate deviations are also discussed for lattice version $X_t^x$. Reference \cite{Landim1992} proves the large deviation principle (LDP) of the occupation time of the SSEP on $\mathbb{Z}^d$ for $d\geq 3$. The $d=2$ case is dealt with in \cite{Lee2004}. Reference \cite{Gao2024} proves a deviation inequality for the SSEP on $\mathbb{Z}^d$. As an application of this deviation inequality, the moderate deviation principle (MDP) of additive functions of SSEP on $\mathbb{Z}^d$ is given, including the MDP of the occupation time as a special case.

Inspired by \cite{Gao2024}, in this paper we give MDPs of the $\mathbb{T}^d$ version centered occupation time $\xi_t^x$. The proof of our main results relies heavily on the duality relationship given in Subsection \ref{subsection 1.2}. For mathematical details, see Section \ref{section two}.

Since 1980s, limit theorems of occupation times are also discussed for other interacting particle systems such as voter models, contact processes, branching Brownian motion, branching random walks, branching $\alpha$-stable processes, binary contact path processes. Readers interested in this topic could resort to References \cite{Birkner2007, Bramson1988, Cox1983, Deuschel1998, Komorowski2012, Li2011, Maillard2009, Schonmann1986, Xue2024}.

\section{Main results}\label{section two}
In this section, we give our main results. From now on, we assume that $p$ is a fixed parameter in $(0,1)$ and $d\geq 2$. We further assume that the SSEP $\{\eta_t\}_{t\geq 0}$ starts from $\nu_p$. Then, for any $x\in \mathbb{T}^d$, the centered occupation time $\xi_t^x$ is defined as in Section \ref{section one}, i.e.,
\[
\xi_t^x=\int_0^t(\eta_s(x)-p)ds.
\]
For later use, we first introduce some notations. For any $x,y\in \mathbb{T}^d$,  we denote by $D(x,y)$ the distance between $x$ and $y$, i.e., $D(x,y)=l$ when and only when there is a self-avoiding path $x=x_0\sim x_1\sim x_2\sim\ldots\sim x_l=y$ on $\mathbb{T}^d$. For any integer $m\geq 1$ and
\[
c=\left(c(1),\ldots,c(m)\right)^T,u=\left(u(1),\ldots,u(d)\right)^T\in \mathbb{R}^m,
\]
where $T$ is the transposition operator, we denote by $c\cdot u$ the inner product of $c$ and $u$, i.e., $c\cdot u=\sum_{i=1}^mc(i)u(i)$. For any integer $m\geq 1$ and $x_1,\ldots,x_m\in \mathbb{T}^d$, we define $\Lambda^{\{x_i\}_{i=1}^m}_t$ as the random vector
\[
\left(\xi_t^{x_1},\ldots,\xi_t^{x_m}\right)^{T}.
\]

To give our main results, we first introduce our rate function. For any integer $m\geq 1$ and $x_1,\ldots,x_m\in \mathbb{T}^d$, we define
\begin{equation}\label{equ rate function}
I_{\{x_i\}_{i=1}^m}(u)=\sup_{c\in \mathbb{R}^m}\left\{c\cdot u-\frac{1}{2}c^T\Gamma_{\{x_i\}_{i=1}^m}c\right\}
\end{equation}
for any $u\in \mathbb{R}^m$, where $\Gamma_{\{x_i\}_{i=1}^m}$ is a $m\times m$ matrix such that
\[
\Gamma_{\{x_i\}_{i=1}^m}(j,k)=2p(1-p)\int_0^{+\infty}p_s(x_j, x_k)ds
\]
for all $1\leq j,k\leq m$, where $\{p_s(\cdot, \cdot)\}_{s\geq 0}$ are the transition probabilities of $\{V_t\}_{t\geq 0}$ defined as in Subsection \ref{subsection 1.2}. Note that, for $d\geq 2$ and $x,y\in \mathbb{T}^d$, $\int_0^{+\infty}p_s(x,y)ds<+\infty$ according to the fact that
\begin{equation}\label{equ 2.1 heat estimation}
p_t(x,y)\leq \left(\sqrt{d}\right)^{-D(x,y)}e^{-t\left(\sqrt{d}-1\right)^2}.
\end{equation}
A proof of Equation \eqref{equ 2.1 heat estimation} is given in Appendix \ref{appendix 1}.

When $m=1$, for any $x\in \mathbb{T}^d$ and $u\in \mathbb{R}$,
\begin{equation}\label{equ 1D rate function}
I_x(u)=\sup_{c\in \mathbb{R}}\left\{cu-\frac{c^2\sigma^2}{2}\right\}=\frac{u^2}{2\sigma^2},
\end{equation}
where
\[
\sigma^2=2p(1-p)\int_0^{+\infty}p_s(x,x)ds.
\]
Note that $\sigma^2$ does not depend on the choice of $x$ according to the spatial homogeneity of our model. Now we give our main result.

\begin{theorem}\label{theorem main mdp}
Let $\{a_t\}_{t\geq 0}$ be a positive sequence such that
\[
\lim_{t\rightarrow+\infty}\frac{a_t}{t}=\lim_{t\rightarrow+\infty}\frac{\sqrt{t}}{a_t}=0.
\]
For any integer $m\geq 1$, $x_1,\ldots,x_m\in \mathbb{T}^d$ and closed set $\mathcal{C}\subseteq \mathbb{R}^m$,
\begin{equation}\label{equ mdp upper bound}
\limsup_{t\rightarrow+\infty}\frac{t}{a_t^2}\log \mathbb{P}_{\nu_p}\left(\frac{1}{a_t}\Lambda_t^{\{x_i\}_{i=1}^m}\in \mathcal{C}\right)\leq -\inf_{u\in \mathcal{C}}I_{\{x_i\}_{i=1}^m}(u).
\end{equation}
For any open set $\mathcal{O}\subseteq \mathbb{R}^m$,
\begin{equation}\label{equ mdp lower bound}
\liminf_{t\rightarrow+\infty}\frac{t}{a_t^2}\log \mathbb{P}_{\nu_p}\left(\frac{1}{a_t}\Lambda_t^{\{x_i\}_{i=1}^m}\in \mathcal{O}\right)\geq -\inf_{u\in \mathcal{O}}I_{\{x_i\}_{i=1}^m}(u).
\end{equation}
\end{theorem}
According to Equation \eqref{equ 1D rate function}, we have the following corollary.
\begin{corollary}\label{corollary 1D mdp}
Let $\{a_t\}_{t\geq 0}$ be a positive sequence such that
\[
\lim_{t\rightarrow+\infty}\frac{a_t}{t}=\lim_{t\rightarrow+\infty}\frac{\sqrt{t}}{a_t}=0.
\]
For any $x\in \mathbb{T}^d$ and $u>0$,
\[
\lim_{t\rightarrow+\infty}\frac{t}{a_t^2}\log \mathbb{P}_{\nu_p}\left(\frac{1}{a_t}\xi_t^x\geq u\right)=-\frac{u^2}{2\sigma^2}.
\]
\end{corollary}

\begin{remark}\label{remark 2.1}
Theorem \ref{theorem main mdp} is consistent with a heuristic covariance analysis, i.e.,
\begin{equation}\label{equ convariance analysis}
\lim_{t\rightarrow+\infty}{\rm Cov}\left(\frac{1}{\sqrt{t}}\xi_t^{x_j}, \frac{1}{\sqrt{t}}\xi_t^{x_k}\right)=\Gamma_{\{x_i\}_{i=1}^m}(j,k).
\end{equation}
Here we give an outline of how to check Equation \eqref{equ convariance analysis}. By Equation \eqref{equ 1.3 1vertexdual}, we have
\[
\mathbb{E}_\eta\eta_t(x)=\sum_{z\in \mathbb{Z}^d}p_t(x,z)\eta(z).
\]
Then, according to the Markov property of $\{\eta_t\}_{t\geq 0}$ and the invariance of $\nu_p$, we have
\[
\lim_{t\rightarrow+\infty}{\rm Cov}\left(\frac{1}{\sqrt{t}}\xi_t^{x_j}, \frac{1}{\sqrt{t}}\xi_t^{x_k}\right)
=2\int_0^{+\infty}\sum_zp_u(x_j,z){\rm Cov}_{\nu_p}\left(\eta(x_k), \eta(z)\right)du.
\]
According to the definition of $\nu_p$, we have ${\rm Cov}_{\nu_p}\left(\eta(x_k), \eta(z)\right)=p(1-p)1_{\{z=x_k\}}$, where $1_{A}$ is the indicator function of event $A$.  Consequently,
\[
\lim_{t\rightarrow+\infty}{\rm Cov}\left(\frac{1}{\sqrt{t}}\xi_t^{x_j}, \frac{1}{\sqrt{t}}\xi_t^{x_k}\right)=2p(1-p)\int_0^{+\infty}p_u(x_j,x_k)du.
\]
\end{remark}

The rest of this paper is organized as follows. In Section \ref{section three}, as preliminaries of the proof of Theorem \ref{theorem main mdp}, we give two replacement lemmas and show that $\{\frac{1}{a_t}\xi_t^x\}_{t\geq 0}$ are exponentially tight. In Sections \ref{section four} and \ref{section five}, we prove Equations \eqref{equ mdp upper bound} and \eqref{equ mdp lower bound} respectively. Our proofs of Equations \eqref{equ mdp upper bound} and \eqref{equ mdp upper bound} utilize the exponential martingale strategy introduced in \cite{Kipnis1989} and the martingale decomposition strategy introduced in \cite{Kipnis1987}, which apply in the SSEP on $\mathbb{T}^d$ due to the replacement lemmas and exponential tightness given in Section \ref{section three}. For mathematical details, see Sections \ref{section three}-\ref{section five}.

\section{Replacement lemmas and exponential tightness}\label{section three}

In this section. we prove two replacement lemmas and the exponential tightness of $\{\frac{1}{a_t}\xi_t^x\}_{t\geq 0}$. To explain the motivation of our replacement lemmas, we first introduce some notations and recall the exponential martingale strategy and the martingale decomposition strategy introduced in \cite{Kipnis1989} and \cite{Kipnis1987} respectively. For any $t\geq 0, \eta\in \{0, 1\}^{\mathbb{T}^d}$ and $x\in \mathbb{T}^d$, we define
\[
G_t^x(\eta)=\sum_{y\in \mathbb{T}^d}(\eta(y)-p)g_t^x(y),
\]
where
\[
g_t^x(y)=\int_0^{+\infty}e^{-\frac{1}{\sqrt{t}}s}p_s(x,y)ds.
\]
Note that $\sum_{y\in \mathbb{T}^d}g_t(y)=\sqrt{t}<+\infty$. By direct calculation, it is easy to check that
\begin{equation}\label{equ 3.1}
\mathcal{L}G_t^x(\eta)=\frac{1}{\sqrt{t}}G_t^x(\eta)-(\eta(x)-p).
\end{equation}
A proof of Equation \eqref{equ 3.1} is given in Appendix \ref{appendix 2}. The martingale decomposition strategy considers the martingale $\{M_s^{t,x}\}_{0\leq s\leq t}$, where
\[
M_s^{t,x}=G_t^x(\eta_s)-G_t^x(\eta_0)-\int_0^s\mathcal{L}G_t^x(\eta_u)du.
\]
By Equation \eqref{equ 3.1}, we have the decomposition
\[
M_s^{t,x}=G_t^x(\eta_s)-G_t^x(\eta_0)-\frac{1}{\sqrt{t}}\int_0^sG_t^x(\eta_u)du+\xi_s^x.
\]
Hence, to utilize the martingale $M_s^{t,x}$ to investigate the limit theorem of $\xi_s^x$, it is inevitable to show that the term $G_t^x(\eta_s)-G_t^x(\eta_0)-\frac{1}{\sqrt{t}}\int_0^sG_t^x(\eta_u)du$ can be neglected after proper scaling, i.e., replaced by $0$. Specific to this paper, since we are concerned with the limit behavior of $\frac{t}{a_t^2}\log\mathbb{P}_{\nu_p}(\frac{1}{a_t}\xi_t^x\in \cdot)$ and $|G_t^x(\eta)|/a_t\leq \frac{\sqrt{t}}{a_t}\rightarrow 0$, it is natural for us to prove the following conclusion, which is our first replacement lemma.
\begin{lemma}\label{lemma 3.1}
For any $\epsilon>0$,
\[
\limsup_{t\rightarrow+\infty}\frac{t}{a_t^2}\log \mathbb{P}_{\nu_p}\Bigg(\Big|\frac{1}{a_t\sqrt{t}}\int_0^tG_t^x(\eta_s)ds\Big|\geq \epsilon\Bigg)=-\infty.
\]
\end{lemma}

In the exponential martingale strategy, it is usual to further consider the martingale $\{\Xi_s^{t,\{x_j\}_{j=1}^m, c}\}_{0\leq s\leq t}$ for $x_1,\ldots,x_m\in \mathbb{T}^d$ and $c=(c_1,\ldots, c_m)^T\in \mathbb{R}^m$, where
\begin{align*}
\Xi_s^{t,\{x_j\}_{j=1}^m, c}=\exp\Bigg\{\frac{a_t}{t}\sum_{j=1}^mc_jG_t^{x_j}(\eta_s)-\frac{a_t}{t}\sum_{j=1}^mc_jG_t^{x_j}(\eta_0)-\int_0^s\frac{\mathcal{L}e^{\frac{a_t}{t}\sum_{j=1}^m
c_jG_t^{x_j}(\eta_u)}}{e^{\frac{a_t}{t}\sum_{j=1}^mc_jG_t^{x_j}(\eta_u)}}du\Bigg\}.
\end{align*}
According to the Taylor's expansion formula up to the second order, it is easy to check that
\begin{align*}
&\Xi_t^{t,\{x_j\}_{j=1}^m, c}\\
&=\exp\Bigg\{\frac{a_t^2}{t}\Big(\frac{1}{a_t}\sum_{j=1}^mc_jM_t^{t,x_j}\\
&\text{\quad}-\frac{1}{2t}\int_0^t\frac{1}{2}\sum_{y\in \mathbb{T}^d}\sum_{z\sim y}\big(\eta_u(z)-\eta_u(y)\big)^2\sum_{j=1}^m\sum_{k=1}^mc_jc_k\left(g_t^{x_j}(y)-g_t^{x_j}(z)\right)\left(g_t^{x_k}(y)-g_t^{x_k}(z)\right)du\\
&\text{\quad}+o(1)\Big)\Bigg\}.
\end{align*}
According to Lemma \ref{lemma 3.1}, $\frac{1}{a_t}M_t^{t,x_j}$ can be replaced by $\frac{1}{a_t}\xi_t^{x_j}$ with a super-exponentially small error. We further require that the term $\big(\eta_u(z)-\eta_u(y)\big)^2$ in the above integral can be replaced by its $\nu_p$-expectation $2p(1-p)$, i.e., we need to prove the following lemma.
\begin{lemma}\label{lemma 3.2}
For any $x,w\in \mathbb{T}^d$, let
\[
\Phi_t^{x,w}=\frac{1}{t}\int_0^t\sum_{y\in \mathbb{T}^d}\sum_{z\sim y}\Big(\big(\eta_u(z)-\eta_u(y)\big)^2-2p(1-p)\Big)\left(g_t^x(y)-g_t^x(z)\right)\left(g_t^w(y)-g_t^w(z)\right)du.
\]
For any $\epsilon>0$ and $x,w\in \mathbb{T}^d$,
\[
\lim_{t\rightarrow+\infty}\frac{t}{a_t^2}\log \mathbb{P}_{\nu_p}\big(|\Phi_t^{x,w}|\geq \epsilon\big)=-\infty.
\]
\end{lemma}
It is easy to check that
\begin{equation}\label{equ 3.2}
\lim_{t\rightarrow+\infty}\sum_{y\in \mathbb{T}^d}\sum_{z\sim y}\left(g_t^x(y)-g_t^x(z)\right)\left(g_t^w(y)-g_t^w(z)\right)=2\int_0^{+\infty}p_s(x,w)ds,
\end{equation}
a proof of which is given in Appendix \ref{appendix 3}. By Lemma \ref{lemma 3.2} and Equation \eqref{equ 3.2}, in the expression of $\Xi_t^{t,\{x_j\}_{j=1}^m, c}$, the term
\[
\frac{1}{2t}\int_0^t\frac{1}{2}\sum_{y\in \mathbb{T}^d}\sum_{z\sim y}\big(\eta_u(z)-\eta_u(y)\big)^2\sum_{j=1}^m\sum_{k=1}^mc_jc_k\left(g_t^{x_j}(y)-g_t^{x_j}(z)\right)\left(g_t^{x_k}(y)-g_t^{x_k}(z)\right)du
\]
can be replaced by $\frac{1}{2}c^T\Gamma_{\{x_i\}_{i=1}^m}c$ with a super-exponentially small error. With above replacement expression of $\Xi_t^{t,\{x_j\}_{j=1}^m, c}$, proofs of Equation \eqref{equ mdp upper bound} for compact sets and Equation \eqref{equ mdp lower bound} for open sets follow from a routine analysis given in literatures such as References \cite{Kipnis1989} and \cite{Gao2003}. To prove Equation \eqref{equ mdp upper bound} for all closed sets, we need to show that $\{\xi_t^x\}_{t\geq 0}$ are exponentially tight, i.e., we require the following lemma.

\begin{lemma}\label{lemma 3.3}
For any $x\in \mathbb{T}^d$,
\begin{equation}\label{equ 3.3}
\limsup_{M\rightarrow+\infty}\limsup_{t\rightarrow+\infty}\frac{t}{a_t^2}\log \mathbb{P}_{\nu_p}\left(\left|\frac{1}{a_t}\xi_t^x\right|\geq M\right)=-\infty.
\end{equation}
\end{lemma}

Our proofs of Lemmas \ref{lemma 3.1}-\ref{lemma 3.3} follow from the same strategy, where Lemma 7.2 in Appendix 1 of \cite{kipnis+landim99}  and duality relationships given in Lemma \ref{lemma dual} play key roles. A detailed proof of Lemma \ref{lemma 3.1} is given in Subsection \ref{subsection 3.1}. Outlines of proofs of Lemmas \ref{lemma 3.2} and \ref{lemma 3.3} are given in Subsection \ref{subsection 3.2}.

\subsection{The proof of Lemma \ref{lemma 3.1}}\label{subsection 3.1}

In this subsection, we prove Lemma \ref{lemma 3.1}. For later use, we first introduce some notations and definitions. For any $t,u\geq 0$, $x\in \mathbb{T}^d$ and $\eta\in \{0, 1\}^{\mathbb{T}^d}$, we denote by $K_t^x(u,\eta)$ the expectation $\mathbb{E}_\eta G_t^x(\eta_u)$. We briefly write $K_t^x(u,\eta)$ as $K_t^x(u)$ when we emphasize that $K_t^x(u,\cdot)$ is a random variable from $\{0,1\}^{\mathbb{T}^d}$ to $\mathbb{R}$. For any $f$ from $\{0, 1\}^{\mathbb{T}^d}$ to $[0, +\infty)$, we call $f$ a $\nu_p$-density if and only if $\int f(\eta)\nu_p(d\eta)=1$. We denote by $\mathfrak{D}$ the Direchlet form of $\{\eta_t\}_{t\geq 0}$, i.e.,
\[
\mathfrak{D}(f)=\frac{1}{4}\int\sum_{x\in \mathbb{T}^d}\sum_{y\sim x}\left(f(\eta^{x,y})-f(\eta)\right)^2\nu_p(d\eta)
\]
for any $f$ from $\{0, 1\}^{\mathbb{T}^d}$ to $\mathbb{R}$. The following lemma plays key role in the proof of Lemma \ref{lemma 3.1}.

\begin{lemma}\label{lemma 3.4}
For any $t>0$ and $x\in \mathbb{T}^d$, under $\nu_p$,
\[
\lim_{u\rightarrow+\infty}K_t^x(u)=0 \text{~in~}L^2.
\]
\end{lemma}

\proof[Proof of Lemma \ref{lemma 3.4}]

According to Lemma \ref{lemma dual}, we have
\begin{equation}\label{equ 3.1.0}
K_t^x(u, \eta)=\sum_{y\in \mathbb{T}^d}\sum_{z\in \mathbb{T}^d}p_u(y,z)(\eta(z)-p)g_t^x(y).
\end{equation}
Hence,
\begin{align*}
\mathbb{E}_{\nu_p}\Big(\big|K_t^x(u)\big|^2\Big)=\sum_{y_1\in \mathbb{T}^d}\sum_{y_2\in \mathbb{T}^d}\sum_{z_1\in \mathbb{T}^d}\sum_{z_2\in \mathbb{T^d}}p_u(y_1, z_1)p_u(y_2, z_2)g_t^x(y_1)g_t^x(y_2)
{\rm Cov}_{\nu_p}\left(\eta(z_1), \eta(z_2)\right).
\end{align*}
Since ${\rm Cov}_{\nu_p}\left(\eta(z_1), \eta(z_2)\right)=p(1-p)1_{\{z_1=z_2\}}$, we have
\begin{align*}
\mathbb{E}_{\nu_p}\Big(\big|K_t^x(u)\big|^2\Big)
&=p(1-p)\sum_{y_1\in \mathbb{T}^d}\sum_{y_2\in \mathbb{T}^d}\sum_{z\in \mathbb{T}^d}p_u(y_1, z)p_u(y_2, z)g_t^x(y_1)g_t^x(y_2)\\
&=p(1-p)\sum_{y_1\in \mathbb{T}^d}\sum_{y_2\in \mathbb{T}^d}p_{2u}(y_1, y_2)g_t^x(y_1)g_t^x(y_2)\\
&=p(1-p)\int_0^{+\infty}\int_0^{+\infty}e^{-\frac{1}{\sqrt{t}}(s_1+s_2)}p_{s_1+s_2+2u}(x,x)ds_1ds_2.
\end{align*}
Consequently, $\lim_{u\rightarrow+\infty}\mathbb{E}_{\nu_p}\Big(\big|K_t^x(u)\big|^2\Big)=0$ according to Equation \eqref{equ 2.1 heat estimation} and the proof is complete.

\qed

Now we prove Lemma \ref{lemma 3.1}.

\proof[Proof of Lemma \ref{lemma 3.1}]

We only show that
\[
\limsup_{t\rightarrow+\infty}\frac{t}{a_t^2}\log \mathbb{P}_{\nu_p}\Bigg(\frac{1}{a_t\sqrt{t}}\int_0^tG_t^x(\eta_s)ds\geq \epsilon\Bigg)=-\infty,
\]
since $\limsup_{t\rightarrow+\infty}\frac{t}{a_t^2}\log \mathbb{P}_{\nu_p}\Bigg(\frac{1}{a_t\sqrt{t}}\int_0^tG_t^x(\eta_s)ds\leq -\epsilon\Bigg)=-\infty$ follows from a similar analysis. According to Markov inequality, for any $\theta>0$,
\[
\mathbb{P}_{\nu_p}\Bigg(\frac{1}{a_t\sqrt{t}}\int_0^tG_t^x(\eta_s)ds\geq \epsilon\Bigg)
\leq e^{-\theta \frac{a_t^2}{t}\epsilon}\mathbb{E}_{\nu_p}\exp\Big\{\frac{a_t\theta}{t^{\frac{3}{2}}}\int_0^tG_t^x(\eta_s)ds\Big\}.
\]
Hence, to complete the proof, we only need to show that, for any $\theta>0$,
\begin{equation}\label{equ 3.1.1}
\limsup_{t\rightarrow+\infty}\frac{t}{a_t^2}\log \mathbb{E}_{\nu_p}\exp\Big\{\frac{a_t\theta}{t^{\frac{3}{2}}}\int_0^tG_t^x(\eta_s)ds\Big\}\leq 0.
\end{equation}
According to Lemma 7.2 in Appendix 1 of \cite{kipnis+landim99},
\[
\mathbb{E}_{\nu_p}\exp\Big\{\frac{a_t\theta}{t^{\frac{3}{2}}}\int_0^tG_t^x(\eta_s)ds\Big\}\leq e^{t\Upsilon_t},
\]
where
\[
\Upsilon_t=\sup_{f \text{~is a~}\nu_p\text{-density}}\Big\{\frac{a_t\theta}{t^{\frac{3}{2}}}\int G_t^x(\eta)f(\eta)\nu_p(d\eta)-\mathfrak{D}(\sqrt{f})\Big\}.
\]
Therefore, to prove Equation \eqref{equ 3.1.1}, we only need to show that, for any $\theta>0$,
\begin{equation}\label{equ 3.1.2}
\limsup_{t\rightarrow+\infty}\sup_{f \text{~is a~}\nu_p\text{-density}}\Big\{\frac{\sqrt{t}\theta}{a_t}\int G_t^x(\eta)f(\eta) \nu_p(d\eta)-\frac{t^2}{a_t^2}\mathfrak{D}(\sqrt{f})\Big\}\leq 0.
\end{equation}
Now we check Equation \eqref{equ 3.1.2}. For any $\nu_p$-density $f$,  according to Lemma \ref{lemma 3.4},
\begin{align*}
\int f(\eta)G_t^x(\eta)\nu_p(d\eta)&=\int f(\eta)K_t^x(0, \eta)\nu_p(d\eta)\\
&=-\int f(\eta)\left(\int_0^{+\infty} \frac{d}{du}K_t^x(u, \eta)du\right)\nu_p(d\eta)\\
&=-\int_0^{+\infty}\left(\int f(\eta)\frac{d}{du}K_t^x(u, \eta)\nu_p(d\eta)\right)du.
\end{align*}
According to the Kolomogrov-Chapman equation, $\frac{d}{du}K_t^x(u, \eta)=\mathcal{L}K_t^x(u, \eta)$. Then, by Equation \eqref{equ 1.2 equilibrium},
\begin{align*}
\int f(\eta)\frac{d}{du}K_t^x(u, \eta)\nu_p(d\eta)&=\int f(\eta)\mathcal{L}K_t^x(u, \eta)\nu_p(d\eta)\\
&=\frac{1}{2}\int K_t^x(u, \eta)\mathcal{L}f(\eta)+f(\eta)\mathcal{L}K_t^x(u, \eta)\nu_p(d\eta)\\
&=\frac{1}{2}\int K_t^x(u, \eta)\mathcal{L}f(\eta)+f(\eta)\mathcal{L}K_t^x(u, \eta)-\mathcal{L}(fK_t^x(u))(\eta)\nu_p(d\eta)\\
&=-\frac{1}{4}\int \sum_{z\in \mathbb{T}^d}\sum_{y\sim z}\left(f(\eta^{z,y})-f(\eta)\right)\left(K_t^x(u, \eta^{z,y})-K_t^x(u, \eta)\right)\nu_p(d\eta).
\end{align*}
Therefore,
\[
\int f(\eta)G_t^x(\eta)\nu_p(d\eta)=\frac{1}{4}\int \sum_{z\in \mathbb{T}^d}\sum_{y\sim z}\Big(f(\eta^{z,y}-f(\eta))\Big)\mathcal{A}_{z,y}^x(\eta)\nu_p(d\eta),
\]
where
\[
\mathcal{A}_{z,y}^x(\eta)=\int_0^{+\infty}K_t^x(u, \eta^{z,y})-K_t^x(u, \eta)du.
\]
Hence, by Cauchy-Schwarz inequality,
\begin{align}\label{equ 3.1.2 and a half}
&\left|\int f(\eta)G_t^x(\eta)\nu_p(d\eta)\right| \notag\\
&=\left|\frac{1}{4}\int \sum_{z\in \mathbb{T}^d}\sum_{y\sim z}\Big(\sqrt{f}(\eta^{z,y})-\sqrt{f}(\eta)\Big)\Big(\sqrt{f}(\eta^{z,y})+\sqrt{f}(\eta)\Big)\mathcal{A}_{z,y}^x(\eta)\nu_p(d\eta)\right|\notag\\
&\leq \frac{1}{4}\sqrt{\int \sum_{z\in \mathbb{T}^d}\sum_{y\sim z}\Big(\sqrt{f}(\eta^{z,y})-\sqrt{f}(\eta)\Big)^2\nu_p(d\eta)}\notag\\
&\text{\quad\quad}\times\sqrt{\int \sum_{z\in \mathbb{T}^d}\sum_{y\sim z}\Big(\sqrt{f}(\eta^{z,y})+\sqrt{f}(\eta)\Big)^2(\mathcal{A}_{z,y}^x(\eta))^2\nu_p(d\eta)}\notag\\
&=\frac{1}{2}\sqrt{\mathfrak{D}(\sqrt{f})}\sqrt{\int \sum_{z\in \mathbb{T}^d}\sum_{y\sim z}\Big(\sqrt{f}(\eta^{z,y})+\sqrt{f}(\eta)\Big)^2(\mathcal{A}_{z,y}^x(\eta))^2\nu_p(d\eta)}.
\end{align}
According to the inequality $(c+b)^2\leq 2c^2+2b^2$,
\begin{align}\label{equ 3.1.3}
&\int \sum_{z\in \mathbb{T}^d}\sum_{y\sim z}\Big(\sqrt{f}(\eta^{z,y})+\sqrt{f}(\eta)\Big)^2(\mathcal{A}_{z,y}^x(\eta))^2\nu_p(d\eta) \notag\\
&\leq 2\int \sum_{z\in \mathbb{T}^d}\sum_{y\sim z}\Big(f(\eta^{z,y})+f(\eta)\Big)(\mathcal{A}_{z,y}^x(\eta))^2\nu_p(d\eta).
\end{align}
By Equation \eqref{equ 3.1.0},
\[
\mathcal{A}_{z,y}^x(\eta)=\int_0^{+\infty}\sum_{w\in \mathbb{T}^d}g_t^x(w)\big(p_u(w,z)-p_u(w,y)\big)\big(\eta(y)-\eta(z)\big)du
\]
and hence
\begin{align*}
(\mathcal{A}_{z,y}^x(\eta))^2
&\leq \int_0^{+\infty}\int_0^{+\infty}\sum_{w_1\in \mathbb{T}^d}\sum_{w_2\in \mathbb{T}^d}\Bigg(g_t^x(w_1)g_t^x(w_2)\\
&\text{\quad\quad}\times\big(p_{u_1}(w_1,z)-p_{u_1}(w_1,y)\big)\big(p_{u_2}(w_2,z)-p_{u_2}(w_2,y)\big)\Bigg)du_1du_2.
\end{align*}
According to an analysis similar with that given in the proof of Equation \eqref{equ 3.2}, we have
\begin{align*}
&\sum_{z\in \mathbb{T}^d}\sum_{y\sim z}\int_0^{+\infty}\int_0^{+\infty}\big(p_{u_1}(w_1,z)-p_{u_1}(w_1,y)\big)\big(p_{u_2}(w_2,z)-p_{u_2}(w_2,y)\big)du_1du_2\\
&=2\sum_{z\in \mathbb{T}^d}\int_0^{+\infty}p_{u_1}(w_1, z)1_{\{w_2=z\}}du_1=2\int_0^{+\infty}p_{u_1}(w_1, w_2)du_1.
\end{align*}
Therefore,
\begin{align}\label{equ 3.1.3 and a half}
\sum_{z\in \mathbb{T}^d}\sum_{y\sim z}(\mathcal{A}_{z,y}^x(\eta))^2&\leq 2\sum_{w_1\in \mathbb{T}^d}\sum_{w_2\in \mathbb{T}^d}g_t^x(w_1)g_t^x(w_2)\int_0^{+\infty}p_u(w_1, w_2)du\notag\\
&=2\int_0^{+\infty}\int_0^{+\infty}\int_0^{+\infty}e^{-\frac{1}{\sqrt{t}}(s_1+s_2)}p_{s_1+s_2+u}(x,x)ds_1ds_2du\notag\\
&\leq  2\int_0^{+\infty}\int_0^{+\infty}\int_0^{+\infty}p_{s_1+s_2+u}(x,x)ds_1ds_2du.
\end{align}
We denote by $J_1$ the term $2\int_0^{+\infty}\int_0^{+\infty}\int_0^{+\infty}p_{s_1+s_2+u}(x,x)ds_1ds_2du$. Note that $J_1<+\infty$ according to Equation \eqref{equ 2.1 heat estimation}. Since $f$ is a $\nu_p$-density,
\[
\int \sum_{z\in \mathbb{T}^d}\sum_{y\sim z}f(\eta)(\mathcal{A}_{z,y}^x(\eta))^2\nu_p(d\eta)\leq J_1.
\]

According to the spatial-homogeneity of $\nu_p$,
\[
\int f(\eta^{y,z})(\mathcal{A}_{z,y}^x(\eta))^2 \nu_p(d\eta)=\int f(\eta)(\mathcal{A}_{z,y}^x(\eta^{y,z}))^2 \nu_p(d\eta)
=\int f(\eta)(\mathcal{A}_{z,y}^x(\eta))^2\nu_p(d\eta).
\]
As a result,
\[
\int \sum_{z\in \mathbb{T}^d}\sum_{y\sim z}\Big(f(\eta^{z,y})+f(\eta)\Big)(\mathcal{A}_{z,y}^x(\eta))^2\nu_p(d\eta)\leq 2J_1.
\]
Then, by Equation \eqref{equ 3.1.3}, we have
\[
\left|\int f(\eta)G_t^x(\eta)\nu_p(d\eta)\right|\leq \sqrt{\mathfrak{D}(\sqrt{f})}\sqrt{J_1}.
\]
Consequently, for any $\nu_p$-density $f$,
\begin{align*}
\frac{\sqrt{t}\theta}{a_t}\int G_t^x(\eta)f(\eta) \nu_p(d\eta)-\frac{t^2}{a_t^2}\mathfrak{D}(\sqrt{f})
&\leq \frac{\sqrt{t}\theta}{a_t}\sqrt{\mathfrak{D}(\sqrt{f})}\sqrt{J_1}-\frac{t^2}{a_t^2}\mathfrak{D}(\sqrt{f})\\
&\leq \sup_{b\in \mathbb{R}}\left\{\frac{\theta}{\sqrt{t}}\sqrt{J_1}b-b^2\right\}=\frac{\theta^2J_1}{4t}.
\end{align*}
Since $\lim_{t\rightarrow+\infty}\frac{\theta^2J_1}{4t}=0$, Equation \eqref{equ 3.1.2} holds and the proof is complete.

\qed

\subsection{Proofs of Lemma \ref{lemma 3.2} and \ref{lemma 3.3}}\label{subsection 3.2}

In this subsection, we prove Lemmas \ref{lemma 3.1} and \ref{lemma 3.2}. To avoid repeating many details similar with those in the proof of Lemma \ref{lemma 3.1}, we only give outlines. The following lemma is an analogue of Lemma \ref{lemma 3.4} and plays key role in the proof of Lemma \ref{lemma 3.2}.

\begin{lemma}\label{lemma 3.2.1}
For any $(x, w)\in \mathbb{Y}, t>0, u\geq 0$ and $\eta\in \{0, 1\}^{\mathbb{T}^d}$, let
\[
\mathcal{U}_t^{x, w}(u, \eta)=\mathbb{E}_{\eta}\sum_{y\in \mathbb{T}^d}\sum_{z\sim y}\Big(\big(\eta_u(z)-\eta_u(y)\big)^2-2p(1-p)\Big)\left(g_t^x(y)-g_t^x(z)\right)\left(g_t^w(y)-g_t^w(z)\right).
\]
Under $\nu_p$,
\[
\lim_{u\rightarrow+\infty}\mathcal{U}_t^{x,w}(u)=0 \text{~in~}L^2.
\]
\end{lemma}

\proof[The outline of the proof of Lemma \ref{lemma 3.2.1}]

For simplicity, we denote by $\hat{g}_{t,y,z}^{x,w}$ the term
\[
\left(g_t^x(y)-g_t^x(z)\right)\left(g_t^w(y)-g_t^w(z)\right).
\]
According to Equation \eqref{equ 1.4 2vertexdual},
\begin{equation}\label{equ mathcalU}
\mathcal{U}_t^{x, w}(u, \eta)=\sum_{y\in \mathbb{T}^d}\sum_{z\sim y}\sum_{(v,r)\in \mathbb{Y}}q_u\Big((y,z), (v,r)\Big)\Big(\left(\eta(v)-\eta(r)\right)^2-2p(1-p)\Big)\hat{g}_{t,y,z}^{x,w}.
\end{equation}
Hence,
\begin{align*}
&\mathbb{E}_{\nu_p}\Big(\left|\mathcal{U}_t^{x,w}(u)\right|^2\Big)\\
&=\sum\mathbb{E}_{\nu_p}\prod_{i=1}^2\Bigg(q_u\Big((y_i,z_i), (v_i,r_i)\Big)\Big(\left(\eta(v_i)-\eta(r_i)\right)^2-2p(1-p)\Big)\hat{g}_{t,y_i,z_i}^{x,w}\Bigg),
\end{align*}
where the sum is over $y_1, y_2\in \mathbb{T}^d, z_1\sim y_1, z_2\sim y_2, (v_1, r_1)\in \mathbb{Y}, (v_2, r_2)\in \mathbb{Y}$. According to the definition of $\nu_p$, we have
\[
\mathbb{E}_{\nu_p}\prod_{i=1}^2\Bigg(q_u\Big((y_i,z_i), (v_i,r_i)\Big)\Big(\left(\eta(v_i)-\eta(r_i)\right)^2-2p(1-p)\Big)\hat{g}_{t,y_i,z_i}^{x,w}\Bigg)=0
\]
when $\left\{v_1, r_1\right\}\bigcap \{v_2, r_2\}=\emptyset$. Hence,
\begin{align*}
\mathbb{E}_{\nu_p}\Big(\left|\mathcal{U}_t^{x,w}(u)\right|^2\Big)
\leq \sum\prod_{i=1}^2\Bigg(q_u\Big((y_i,z_i), (v_i,r_i)\Big)\hat{g}_{t,y_i,z_i}^{x,w}\Bigg),
\end{align*}
where the sum is over $y_1, y_2\in \mathbb{T}^d, z_1\sim y_1, z_2\sim y_2$ and $(v_1, r_1)\in \mathbb{Y}, (v_2, r_2)\in \mathbb{Y}$ such that
\[
\left\{v_1, r_1\right\}\bigcap \{v_2, r_2\}\neq \emptyset.
\]
According to the definition of $q_u\left(\cdot, \cdot\right)$, $\sum_{(v_1, r_1), (v_2, r_2)\in \mathbb{Y}\atop \left\{v_1, r_1\right\}\bigcap \{v_2, r_2\}\neq \emptyset}\prod_{i=1}^2q_u\Big((y_i,z_i), (v_i,r_i)\Big)$ is the probability of the event
\[
\left\{Y_u^{y_1, z_1}(1), Y_u^{y_1, z_1}(2)\right\}\bigcap \left\{\hat{Y}_u^{y_2, z_2}(1), \hat{Y}_u^{y_2, z_2}(2)\right\}\neq \emptyset,
\]
where $\{\hat{Y}_t\}_{t\geq 0}$ is an independent copy of $\{Y_t\}_{t\geq 0}$. According the definition of $\{Y_t\}_{t\geq 0}$, $\{Y_t(k)\}_{t\geq 0}$ is a copy of $\{V_t\}_{t\geq 0}$ for $k=1,2$. Hence, for $k=1,2$ and $l=1,2$,
\begin{align*}
\mathbb{P}\left(Y_u^{y_1, z_1}(l)=\hat{Y}^{y_2, z_2}_u(k)\right)&=\sum_{r\in \mathbb{T}^d}p_u\Big((y_1, z_1)(l),r\Big)p_u\Big((y_2, z_2)(k),r\Big)\\
&=p_{2u}\Big((y_1, z_1)(l), (y_2, z_2)(k)\Big)\leq e^{-2u\left(\sqrt{d}-1\right)^2}
\end{align*}
according to Equation \eqref{equ 2.1 heat estimation}. Therefore,
\[
\sum_{(v_1, r_1), (v_2, r_2)\in \mathbb{Y}\atop \left\{v_1, r_1\right\}\bigcap \{v_2, r_2\}\neq \emptyset}\prod_{i=1}^2q_u\Big((y_i,z_i), (v_i,r_i)\Big)
\leq 4e^{-2u\left(\sqrt{d}-1\right)^2}
\]
and hence
\[
\mathbb{E}_{\nu_p}\Big(\left|\mathcal{U}_t^{x,w}(u)\right|^2\Big)\leq 4e^{-2u\left(\sqrt{d}-1\right)^2}\Big(\sum_{y\in \mathbb{T}^d}\sum_{z\sim y}\hat{g}_{t,y,z}^{x,w}\Big)^2.
\]
According to the proof of Equation \eqref{equ 3.2},
\begin{align*}
\Big(\sum_{y\in \mathbb{T}^d}\sum_{z\sim y}\hat{g}_{t,y,z}^{x,w}\Big)^2=\Bigg(2g_t^x(w)-\frac{2}{\sqrt{t}}\int_0^{+\infty}\int_0^{+\infty}e^{-\frac{1}{\sqrt{t}}(s_1+s_2)}p_{s_1+s_2}(x,x)ds_1ds_2\Bigg)^2<+\infty
\end{align*}
and hence Lemma \ref{lemma 3.2.1} holds.

\qed

Now we give the proof of Lemma \ref{lemma 3.2}.

\proof[The outline of the proof of Lemma \ref{lemma 3.2}]

According to Markov inequality and Lemma 7.2 in Appendix 1 of \cite{kipnis+landim99}, we only need to show that, for any $\theta>0$,
\begin{equation}\label{equ 3.2.101}
\limsup_{t\rightarrow+\infty}\sup_{f \text{~is a~}\nu_p\text{-density}}\Big\{\theta\int \mathcal{U}_t^{x,w}(0, \eta)f(\eta) \nu_p(d\eta)-\frac{t^2}{a_t^2}\mathfrak{D}(\sqrt{f})\Big\}\leq 0.
\end{equation}
According to Lemma \ref{lemma 3.2.1} and an analysis similar with that leading to Equation \eqref{equ 3.1.2 and a half}, we have
\begin{align*}
&\left|\int f(\eta)\mathcal{U}_t^{x,w}(0, \eta)\nu_p(d\eta)\right| \\
&\leq \frac{1}{2}\sqrt{\mathfrak{D}(\sqrt{f})}\sqrt{\int \sum_{z\in \mathbb{T}^d}\sum_{y\sim z}\Big(\sqrt{f}(\eta^{z,y})+\sqrt{f}(\eta)\Big)^2(\mathcal{H}_{z,y}^x(\eta))^2\nu_p(d\eta)},
\end{align*}
where
\[
\mathcal{H}_{z,y}^x(\eta)=\int_0^{+\infty}\mathcal{U}_t^{x,w}(u, \eta^{z,y})-\mathcal{U}_t^{x,w}(u, \eta)du.
\]
By Equation \eqref{equ mathcalU},
\[
\left|\mathcal{H}_{z,y}^x(\eta)\right|\leq \int_0^{+\infty}\sum q_u\Big((\omega, \tau), (v,r)\Big)\hat{g}_{t,\omega,\tau}^{x,w}du,
\]
where the sum is over $\omega\in \mathbb{T}^d, \tau\sim \omega$ and $(v,r)\in \mathbb{Y}$ such that $\{v,r\}\bigcap\{y,z\}\neq \emptyset$. Therefore,
\begin{align*}
\sum_{y\in \mathbb{T}^d}\sum_{z\sim y}\left|\mathcal{H}_{z,y}^x(\eta)\right|^2\leq \int_0^{+\infty}\int_0^{+\infty}\sum \prod_{i=1}^2\Bigg(q_{u_i}\Big((\omega_i, \tau_i), (v_i,r_i)\Big)\hat{g}_{t,\omega_i,\tau_i}^{x,w}\Bigg)du_1du_2,
\end{align*}
where the sum in the integral is over $y\in \mathbb{T}^d, z\sim y$, $\omega_1, \omega_2\in \mathbb{T}^d$, $\tau_1\sim \omega_1, \tau_2\sim \omega_2$ and
\[
(v_1, r_1), (v_2, r_2)\in \mathbb{Y}
\]
such that
$\{v_i,r_i\}\bigcap\{y,z\}\neq \emptyset$ for $i=1,2$. According to the definition of $q_u(\cdot, \cdot)$, the sum
\[
\sum_{y\in \mathbb{T}^d}\sum_{z\sim y}\sum_{(v_1, r_1), (v_2, r_2)\in \mathbb{Y}, \atop \{v_i, r_i\}\bigcap \{y,z\}\neq \emptyset \text{~for~}i=1,2}\prod_{i=1}^2q_{u_i}\Big((\omega_i, \tau_i), (v_i,r_i)\Big)
\]
is at most the probability of the event
\[
D\Big(Y_{u_1}^{w_1, \tau_1}(k), \hat{Y}_{u_2}^{w_2, \tau_2}(l)\Big)\leq 1
\]
for some $l=1,2$ and $k=1,2$, where $\{\hat{Y}_t\}_{t\geq 0}$ is an independent copy of $\{Y_t\}_{t\geq 0}$. Then, according to Equation \eqref{equ 2.1 heat estimation}, the total probability formula and the fact that each vertex on $\mathbb{T}^d$ has $d+1$ neighbors, we have
\begin{align*}
&\sum_{y\in \mathbb{T}^d}\sum_{z\sim y}\sum_{(v_1, r_1), (v_2, r_2)\in \mathbb{Y}, \atop \{v_i, r_i\}\bigcap \{y,z\}\neq \emptyset \text{~for~}i=1,2}\prod_{i=1}^2q_{u_i}\Big((\omega_i, \tau_i), (v_i,r_i)\Big)\\
&\leq 4(d+2)\min\{e^{-u_1\left(\sqrt{d}-1\right)^2}, e^{-u_2\left(\sqrt{d}-1\right)^2}\}.
\end{align*}
Hence, according to the proof of Equation of \eqref{equ 3.2}, we have
\[
\sum_{y\in \mathbb{T}^d}\sum_{z\sim y}\left|\mathcal{H}_{z,y}^x(\eta)\right|^2\leq J_3,
\]
where
\[
J_3=\left(2\int_0^{+\infty}p_s(x,w)ds\right)^2\int_0^{+\infty}\int_0^{+\infty} 4(d+2)e^{-\max\{u_1, u_2\}\left(\sqrt{d}-1\right)^2}du_1du_2<+\infty.
\]
Consequently, for any $\nu_p$-density $f$,
\begin{align*}
\theta\int \mathcal{U}_t^{x,w}(0, \eta)f(\eta) \nu_p(d\eta)-\frac{t^2}{a_t^2}\mathfrak{D}(\sqrt{f})
&\leq \theta\sqrt{\mathfrak{D}(\sqrt{f})}\sqrt{J_3}-\frac{t^2}{a_t^2}\mathfrak{D}(\sqrt{f})\\
&\leq \sup_{b\in \mathbb{R}}\Big(\theta\sqrt{J_3}b-\frac{t^2}{a_t^2}b^2\Big)=\frac{\theta^2a_t^2J_3}{4t^2}.
\end{align*}
Since $a_t/t\rightarrow 0$, Equation \eqref{equ 3.2.101} holds and the proof is complete.

\qed

Now we prove Lemma \ref{lemma 3.3}. We need the following lemma as a preliminary, which is an analogue of Lemmas \ref{lemma 3.4} and .

\begin{lemma}\label{lemma 3.2.2}
For any $x\in \mathbb{T}^d, t>0, u\geq 0$ and $\eta\in \{0, 1\}^{\mathbb{T}^d}$, let $\mathcal{R}_t^x(u, \eta)=\mathbb{E}_{\eta}\left(\eta_u(x)-p\right)$.
Under $\nu_p$,
\[
\lim_{u\rightarrow+\infty}\mathcal{R}_t^x(u)=0 \text{~in~}L^2.
\]
\end{lemma}

\proof[The outline of the proof of Lemma \ref{lemma 3.2.2}]

According to Equation  \eqref{equ 1.3 1vertexdual},
\[
\mathcal{R}_t^x(u, \eta)=\sum_{y\in \mathbb{T}^d}p_u(x,y)(\eta(y)-p).
\]
Hence,
\begin{align*}
\mathbb{E}_{\nu_p}\Big(\big|\mathcal{R}_t^x(u)\big|^2\Big)&=\sum_{y_1\in \mathbb{T}^d}\sum_{y_2\in \mathbb{T}^d}p_u(x, y_1)p_u(x, y_2){\rm Cov}_{\nu_p}\left(\eta(y_1), \eta(y_2)\right)\\
&=p(1-p)\sum_{y\in \mathbb{T}^d}p_u(x, y)p_u(x, y)=p(1-p)p_{2u}(x, x).
\end{align*}
Therefore, Lemma \ref{lemma 3.2.2} follows from Equation \eqref{equ 2.1 heat estimation}.

\qed

At last, we give the proof of Lemma \ref{lemma 3.3}.

\proof[The outline of the proof of Lemma \ref{lemma 3.3}]

According to Markov inequality and Lemma 7.2 in Appendix 1 of \cite{kipnis+landim99}, we only need to show that
\begin{equation}\label{equ 3.2.201}
\limsup_{t\rightarrow+\infty}\sup_{f \text{~is a~}\nu_p\text{-density}}\Big\{\frac{t}{a_t}\int (\eta(x)-p)f(\eta) \nu_p(d\eta)-\frac{t^2}{a_t^2}\mathfrak{D}(\sqrt{f})\Big\}<+\infty.
\end{equation}
According to Lemma \ref{lemma 3.2.2} an analysis similar with that leading to Equation \eqref{equ 3.1.2 and a half}, we have
\begin{align*}
&\left|\int f(\eta)(\eta(x)-p)\nu_p(d\eta)\right| \\
&\leq \frac{1}{2}\sqrt{\mathfrak{D}(\sqrt{f})}\sqrt{\int \sum_{z\in \mathbb{T}^d}\sum_{y\sim z}\Big(\sqrt{f}(\eta^{z,y})+\sqrt{f}(\eta)\Big)^2(\mathcal{B}_{z,y}^x(\eta))^2\nu_p(d\eta)},
\end{align*}
where
\[
\mathcal{B}_{z,y}^x(\eta)=\int_0^{+\infty}\mathcal{R}_t^x(u, \eta^{z,y})-\mathcal{R}_t^x(u, \eta)du.
\]
According to an analysis similar with that leading to Equation \eqref{equ 3.1.3 and a half}, we have
\[
\sum_{z\in \mathbb{T}^d}\sum_{y\sim z}(\mathcal{B}_{z,y}^x(\eta))^2\leq 2\int_0^{+\infty}p_u(x,x)du.
\]
Consequently,
\[
\left|\int f(\eta)(\eta(x)-p)\nu_p(d\eta)\right|\leq \sqrt{\mathfrak{D}(\sqrt{f})}\sqrt{J_2},
\]
where $J_2=2\int_0^{+\infty}p_u(x,x)du$. Hence, for any $\nu_p$-density $f$,
\begin{align*}
\frac{t}{a_t}\int (\eta(x)-p)f(\eta) \nu_p(d\eta)-\frac{t^2}{a_t^2}\mathfrak{D}(\sqrt{f})& \leq \frac{t}{a_t}\sqrt{\mathfrak{D}(\sqrt{f})}\sqrt{J_2}-\frac{t^2}{a_t^2}\mathfrak{D}(\sqrt{f})\\
&\leq \sup_{b\in \mathbb{R}}\Big(b\sqrt{J_2}-b^2\Big)=\frac{J_2}{4}.
\end{align*}
As a result, Equation \eqref{equ 3.2.201} holds and the proof is complete.

\qed

\section{The proof of Equation \eqref{equ mdp upper bound}}\label{section four}

Now we give the proof of Equation \eqref{equ mdp upper bound}.

\proof[Proof of Equation \eqref{equ mdp upper bound}]
Throughout this proof we assume that $x_1, x_2, \ldots, x_m$ are fixed.
According to Lemma \ref{lemma 3.3}, we only need to prove Equation \eqref{equ mdp upper bound} for all compact $\mathcal{C}\subseteq \mathbb{R}^m$. For any $c=(c_1,\ldots, c_m)^T\in \mathbb{R}^m$, we define
\[
\Xi_s^{t,\{x_j\}_{j=1}^m, c}=\exp\Bigg\{\frac{a_t}{t}\sum_{j=1}^mc_jG_t^{x_j}(\eta_s)-\frac{a_t}{t}\sum_{j=1}^mc_jG_t^{x_j}(\eta_0)-\int_0^s\frac{\mathcal{L}e^{\frac{a_t}{t}\sum_{j=1}^m
c_jG_t^{x_j}(\eta_u)}}{e^{\frac{a_t}{t}\sum_{j=1}^mc_jG_t^{x_j}(\eta_u)}}du\Bigg\}
\]
as in Section \ref{section three}. According to Feynman-Kac formula, $\{\Xi_s^{t, \{x_j\}_{j=1}^m, c}\}_{0\leq s\leq t}$ is a martingale. According to the definition of $\mathcal{L}$, we have
\begin{align*}
\frac{\mathcal{L}e^{\frac{a_t}{t}\sum_{j=1}^m
c_jG_t^{x_j}(\eta_u)}}{e^{\frac{a_t}{t}\sum_{j=1}^mc_jG_t^{x_j}(\eta_u)}}=\frac{1}{2}\sum_{y\in \mathbb{T}^d}\sum_{z\sim y}\Bigg(e^{\frac{a_t}{t}(\eta_u(z)-\eta_u(y))\sum_{j=1}^mc_j\big(g_t^{x_j}(y)-g_t^{x_j}(z)\big)}-1\Bigg).
\end{align*}
Then, according to the Taylor's expansion formula up to the second order,
\[
\int_0^t\frac{\mathcal{L}e^{\frac{a_t}{t}\sum_{j=1}^m
c_jG_t^{x_j}(\eta_u)}}{e^{\frac{a_t}{t}\sum_{j=1}^mc_jG_t^{x_j}(\eta_u)}}du
=\frac{a_t^2}{t}\left({\rm \uppercase\expandafter{\romannumeral1}}+{\rm \uppercase\expandafter{\romannumeral2}}+o(1)\right),
\]
where
\begin{align*}
{\rm \uppercase\expandafter{\romannumeral1}}&=\frac{1}{a_t}\int_0^t\sum_{j=1}^mc_j\frac{1}{2}\sum_{y\in \mathbb{T}^d, z\sim y}(\eta_u(z)-\eta_u(y))
(g_t^{x_j}(y)-g_t^{x_j}(z))du\\
&=\frac{1}{a_t}\int_0^t\sum_{j=1}^mc_j\mathcal{L}G_t^{x_j}(\eta_u)du
=\frac{1}{a_t\sqrt{t}}\sum_{j=1}^mc_j\int_0^tG_t^{x_j}(\eta_u)du-\sum_{j=1}^mc_j\frac{1}{a_t}\xi_t^{x_j}
\end{align*}
according to Equation \eqref{equ 3.1} and
\begin{align*}
{\rm \uppercase\expandafter{\romannumeral2}}
=\frac{1}{2t}\int_0^t\frac{1}{2}\sum_{y\in \mathbb{T}^d}\sum_{z\sim y}\big(\eta_u(z)-\eta_u(y)\big)^2\sum_{j=1}^m\sum_{k=1}^mc_jc_k\left(g_t^{x_j}(y)-g_t^{x_j}(z)\right)\left(g_t^{x_k}(y)-g_t^{x_k}(z)\right)du.
\end{align*}
Therefore, according to Equation \eqref{equ 3.2} and Lemmas \ref{lemma 3.1}, \ref{lemma 3.2}, we have
\[
{\rm \uppercase\expandafter{\romannumeral1}}+{\rm \uppercase\expandafter{\romannumeral2}}
=-\sum_{j=1}^mc_j\frac{1}{a_t}\xi_t^{x_j}+\frac{1}{2}c^T\Gamma_{\{x_j\}_{j=1}^m}c+\varepsilon_{1,t},
\]
where $\limsup_{t\rightarrow+\infty}\frac{t}{a_t^2}\log \mathbb{P}_{\nu_p}\left(|\varepsilon_{1,t}|\geq \epsilon\right)=-\infty$ for any $\epsilon>0$. Then, since $\frac{|G_t^{x_j}(\eta)|}{a_t}\leq \frac{\sqrt{t}}{a_t}$, we have
\begin{equation}\label{equ 4.1}
\Xi_t^{t,\{x_j\}_{j=1}^m, c}=\exp\left\{\frac{a_t^2}{t}\Big(\sum_{j=1}^mc_j\frac{1}{a_t}\xi_t^{x_j}-\frac{1}{2}c^T\Gamma_{\{x_j\}_{j=1}^m}c+\varepsilon_{2,t}\Big)\right\},
\end{equation}
where $\limsup_{t\rightarrow+\infty}\frac{t}{a_t^2}\log \mathbb{P}_{\nu_p}\left(|\varepsilon_{2,t}|\geq \epsilon\right)=-\infty$ for any $\epsilon>0$. Therefore, for any $\epsilon>0$ and compact $\mathcal{C}\subseteq \mathbb{R}^m$,
\begin{align*}
1&=\mathbb{E}_{\nu_p}\Xi_t^{t,\{x_j\}_{j=1}^m, c}\geq \mathbb{E}_{\nu_p}\Bigg(\Xi_t^{t,\{x_j\}_{j=1}^m, c}1_{\{\frac{1}{a_t}\Lambda_t^{\{x_i\}_{i=1}^m}\in \mathcal{C}, |\epsilon_{2,t}|\leq \epsilon\}}\Bigg)\\
&\geq \mathbb{P}_{\nu_p}\left(\frac{1}{a_t}\Lambda_t^{\{x_i\}_{i=1}^m}\in \mathcal{C}, |\epsilon_{2,t}|\leq \epsilon\right)
e^{\frac{a_t^2}{t}\inf_{u\in \mathcal{C}}\{c\cdot u-\frac{1}{2}c^T\Gamma_{\{x_j\}_{j=1}^m}c-\epsilon\}}.
\end{align*}
Then, since $\limsup_{t\rightarrow+\infty}\frac{t}{a_t^2}\log \mathbb{P}_{\nu_p}\left(|\varepsilon_{2,t}|\geq \epsilon\right)=-\infty$, we have
\begin{align*}
&\limsup_{t\rightarrow+\infty}\frac{t}{a_t^2}\log \mathbb{P}_{\nu_p}\left(\frac{1}{a_t}\Lambda_t^{\{x_i\}_{i=1}^m}\in \mathcal{C}\right)=\limsup_{t\rightarrow+\infty}\frac{t}{a_t^2}\log \mathbb{P}_{\nu_p}\left(\frac{1}{a_t}\Lambda_t^{\{x_i\}_{i=1}^m}\in \mathcal{C}, |\epsilon_{2,t}|\leq \epsilon\right)\\
&\leq -\inf_{u\in \mathcal{C}}\{c\cdot u-\frac{1}{2}c^T\Gamma_{\{x_j\}_{j=1}^m}c-\epsilon\}.
\end{align*}
Since $\epsilon$ and $c$ are arbitrary, we have
\[
\limsup_{t\rightarrow+\infty}\frac{t}{a_t^2}\log \mathbb{P}_{\nu_p}\left(\frac{1}{a_t}\Lambda_t^{\{x_i\}_{i=1}^m}\in \mathcal{C}\right)
\leq -\sup_{c\in \mathbb{R}^m}\inf_{u\in \mathcal{C}}\{c\cdot u-\frac{1}{2}c^T\Gamma_{\{x_j\}_{j=1}^m}c\}.
\]
Since $c\cdot u-\frac{1}{2}c^T\Gamma_{\{x_j\}_{j=1}^m}c$ is linear in $u$ and concave in $c$, according to the minimax theorem given in \cite{Sion1958}, we have
\[
\sup_{c\in \mathbb{R}^m}\inf_{u\in \mathcal{C}}\{c\cdot u-\frac{1}{2}c^T\Gamma_{\{x_j\}_{j=1}^m}c\}
=\inf_{u\in \mathcal{C}}\sup_{c\in \mathbb{R}^m}\{c\cdot u-\frac{1}{2}c^T\Gamma_{\{x_j\}_{j=1}^m}c\}=\inf_{u\in \mathcal{C}}I_{\{x_i\}_{i=1}^m}(u)
\]
and the proof is complete.

\qed

\section{The proof of Equation \eqref{equ mdp lower bound}}\label{section five}

In this section, we prove Equation \eqref{equ mdp lower bound}. We first give two lemmas as preliminaries.

\begin{lemma}\label{lemma 5.1}
If $u\in \mathbb{R}^m$ makes $I_{\{x_i\}_{i=1}^m}(u)<+\infty$, then there exists $\varphi\in \mathbb{R}^m$ such that
$\Gamma_{\{x_j\}_{j=1}^m}\varphi=u$ and
\[
I_{\{x_i\}_{i=1}^m}(u)=\varphi\cdot u-\frac{1}{2}\varphi^T\Gamma_{\{x_j\}_{j=1}^m}\varphi=\frac{1}{2}\varphi^T\Gamma_{\{x_j\}_{j=1}^m}\varphi.
\]
\end{lemma}

The proof of Lemma \ref{lemma 5.1} follows from a routine analysis utilizing Riesz representation theorem, the detail of which we omit in this paper.

\begin{lemma}\label{lemma 5.2}
For given $c\in \mathbb{R}^m$ and $x_1,\ldots, x_m\in \mathbb{T}^d$, let $\hat{\mathbb{P}}_{c,t}^{\{x_i\}_{i=1}^m}$ be the probability measure such that
\[
\frac{d\hat{\mathbb{P}}_{c,t}^{\{x_i\}_{i=1}^m}}{d\mathbb{P}_{\nu_p}}=\Xi_t^{t,\{x_j\}_{j=1}^m, c}.
\]
The process $\{\frac{1}{a_t}\Lambda_t^{\{x_i\}_{i=1}^m}\}_{t\geq 0}$ converges in $\hat{\mathbb{P}}_{c,t}^{\{x_i\}_{i=1}^m}$-probability to $\Gamma_{\{x_j\}_{j=1}^m}c$ as $t\rightarrow+\infty$.
\end{lemma}

Lemma \ref{lemma 5.2} is a standard step in the proof of the MDP lower bound, the proof of which is given in Appendix \ref{appendix 4}. At last, we give the proof of Equation \eqref{equ mdp lower bound}.

\proof[Proof of Equation \eqref{equ mdp lower bound}]

Equation \eqref{equ mdp lower bound} is trivial when $\inf_{u\in \mathcal{O}}I_{\{x_i\}_{i=1}^m}(u)=+\infty$. Hence, we only deal with the case where $\inf_{u\in \mathcal{O}}I_{\{x_i\}_{i=1}^m}(u)<+\infty$. For any $\epsilon>0$, there exists $u_\epsilon\in \mathcal{O}$ such that
\[
I_{\{x_i\}_{i=1}^m}(u_\epsilon)\leq \inf_{u\in \mathcal{O}}I_{\{x_i\}_{i=1}^m}(u)+\epsilon.
\]
By Lemma \ref{lemma 5.1}, there exists $\varphi_\epsilon\in \mathbb{R}^m$ such that $\Gamma_{\{x_j\}_{j=1}^m}\varphi_\epsilon=u_\epsilon$ and
\[
I_{\{x_i\}_{i=1}^m}(u_\epsilon)=\varphi_\epsilon\cdot u_\epsilon-\frac{1}{2}\varphi_\epsilon^T\Gamma_{\{x_j\}_{j=1}^m}\varphi_\epsilon=\frac{1}{2}\varphi_\epsilon^T\Gamma_{\{x_j\}_{j=1}^m}\varphi_\epsilon.
\]
We define $D_1$ as
\[
D_1=\left\{u\in \mathbb{R}^m:~\Big|\varphi_\epsilon\cdot u-\frac{1}{2}\varphi_\epsilon^T\Gamma_{\{x_j\}_{j=1}^m}\varphi_\epsilon-I_{\{x_i\}_{i=1}^m}(u_\epsilon)\Big|<\epsilon\right\}
\]
and denote by $E_{1,t}$ the event $\left\{\frac{1}{a_t}\Lambda_t^{\{x_i\}_{i=1}^m}\in D_1\bigcap \mathcal{O}\right\}$. By Lemma \ref{lemma 5.2}, $\frac{1}{a_t}\Lambda_t^{\{x_i\}_{i=1}^m}$ converges in $\hat{\mathbb{P}}_{\varphi_\epsilon,t}^{\{x_i\}_{i=1}^m}$-probability to $u_\epsilon$ and hence
\[
\lim_{t\rightarrow+\infty}\hat{\mathbb{P}}_{\varphi_{\epsilon},t}^{\{x_i\}_{i=1}^m}\Big(E_{1,t}\Big)=1.
\]
By Taylor's expansion formula up to the second order and Equation \eqref{equ 3.2}, there exists a constant $M<+\infty$ independent of $t$ such that
\[
\left(\Xi_t^{t,\{x_j\}_{j=1}^m, \varphi_\epsilon}\right)^2\leq \Xi_t^{t,\{x_j\}_{j=1}^m, 2\varphi_\epsilon}\exp\left\{\frac{a_t^2}{t}M\right\}
\]
for sufficiently large $t$. Hence, by Cauchy-Schwarz inequality,
\[
\limsup_{t\rightarrow+\infty}\frac{t}{a_t^2}\log \hat{\mathbb{P}}_{\varphi_{\epsilon},t}^{\{x_i\}_{i=1}^m}\left(|\varepsilon_{2,t}|\geq \epsilon\right)=-\infty,
\]
where $\varepsilon_{2,t}$ is defined as in Equation \eqref{equ 4.1}. Hence,
\[
\lim_{t\rightarrow+\infty}\hat{\mathbb{P}}_{\varphi_{\epsilon},t}^{\{x_i\}_{i=1}^m}\Big(E_{2,t}\Big)=1,
\]
where $E_{2,t}=E_{1,t}\bigcap \left\{|\varepsilon_{2,t}|<\epsilon\right\}$. Consequently, by Equation \eqref{equ 4.1},
\begin{align*}
\mathbb{P}_{\nu_p}\left(\frac{1}{a_t}\Lambda_t^{\{x_i\}_{i=1}^m}\in \mathcal{O}\right)
&=\hat{\mathbb{E}}_{\varphi_{\epsilon},t}^{\{x_i\}_{i=1}^m}\Bigg(\Big(\Xi_t^{t,\{x_j\}_{j=1}^m, \varphi_\epsilon}\Big)^{-1}1_{\{\frac{1}{a_t}\Lambda_t^{\{x_i\}_{i=1}^m}\in \mathcal{O}\}}\Bigg)\\
&\geq \hat{\mathbb{E}}_{\varphi_{\epsilon},t}^{\{x_i\}_{i=1}^m}\Bigg(\Big(\Xi_t^{t,\{x_j\}_{j=1}^m, \varphi_\epsilon}\Big)^{-1}1_{E_{2,t}}\Bigg)\\
&\geq \exp\left\{-\frac{a_t^2}{t}\left(I_{\{x_i\}_{i=1}^m}(u_\epsilon)+2\epsilon\right)\right\}
\hat{\mathbb{P}}_{\varphi_{\epsilon},t}^{\{x_i\}_{i=1}^m}\Big(E_{2,t}\Big)\\
&=\exp\left\{-\frac{a_t^2}{t}\left(I_{\{x_i\}_{i=1}^m}(u_\epsilon)+2\epsilon\right)\right\}(1+o(1)).
\end{align*}
Therefore,
\begin{align*}
\liminf_{t\rightarrow+\infty}\frac{t}{a_t^2}\log \mathbb{P}_{\nu_p}\left(\frac{1}{a_t}\Lambda_t^{\{x_i\}_{i=1}^m}\in \mathcal{O}\right)
&\geq -I_{\{x_i\}_{i=1}^m}(u_\epsilon)-2\epsilon\\
&\geq -\inf_{u\in \mathcal{O}}I_{\{x_i\}_{i=1}^m}(u)-3\epsilon.
\end{align*}
Since $\epsilon$ is arbitrary, the proof is complete.

\qed

\quad

\appendix{}
\section{Appendix}
\subsection{The proof of Equation \eqref{equ 2.1 heat estimation}}\label{appendix 1}

\proof[Proof of Equation \eqref{equ 2.1 heat estimation}]

Suppose that $D(x,y)=k$. According to the structure of $\mathbb{T}^d$, there is a function $\beta$ from $\mathbb{T}^d$ to $\mathbb{Z}$ such that, for each $z\in \mathbb{T}^d$, one neighbor $w$ of $z$ satisfies $\beta(w)=\beta(z)-1$ and other $d$ neighbors $v$ of $z$ satisfies $\beta(v)=\beta(z)+1$. Without loss of generality, we assume that $\beta(x)=0$ and $\beta(y)=k$. Then, according to the spatial homogeneity of $\{V_t\}_{t\geq 0}$,
\begin{equation}\label{equ A.1}
p_t(x,y)=\frac{1}{d^k}\mathbb{P}\left(\beta(V_t^x)=k, D(V_t^x, x)=k\right)\leq \frac{1}{d^k}\mathbb{P}\left(\beta(V_t^x)=k\right).
\end{equation}
Since $\{\beta(V_t)\}_{t\geq 0}$ is a continuous-time Markov process such that $\beta(V_t)\rightarrow \beta(V_t)+1$ at rate $d$ and $\beta(V_t)\rightarrow \beta(V_t)-1$ at rate $1$, we have
\[
\frac{d}{dt}\mathbb{E}e^{-\theta \beta(V_t^x)}=\Big(d(e^{-\theta}-1)+(e^\theta-1)\Big)\mathbb{E}e^{-\theta \beta(V_t^x)}
\]
for any $\theta>0$ according to Kolmogorov-Chapman equation. Hence,
\[
\mathbb{E}e^{-\theta \beta(V_t^x)}=e^{t\big(d(e^{-\theta}-1)+(e^\theta-1)\big)}.
\]
By Markov inequality,
\[
\mathbb{P}\left(\beta(V_t^x)=k\right)\leq \mathbb{P}\left(\beta(V_t^x)\leq k\right)\leq e^{\theta k}\mathbb{E}e^{-\theta \beta(V_t^x)}.
\]
Taking $\theta=\log (\sqrt{d})$, then Equation \eqref{equ 2.1 heat estimation} follows from Equation \eqref{equ A.1}.

\qed

\subsection{The proof of Equation \eqref{equ 3.1}}\label{appendix 2}

\proof[Proof of Equation \eqref{equ 3.1}]

According to the definition of $\mathcal{L}$,
\begin{align*}
\mathcal{L}G_t^x(\eta)&=\sum_{y\in \mathbb{T}^d}\mathcal{L}(\eta(y)-p)g_t^x(y)=\sum_{y\in \mathbb{T}^d}\sum_{z\sim y}\Big((\eta(z)-p)-(\eta(y)-p)\Big)g_t^x(y)\\
&=\sum_{y\in \mathbb{T}^d}(\eta(y)-p)\Big(\sum_{z\sim y}g_t^x(z)-(d+1)g_t^x(y)\Big).
\end{align*}
For any $y\in \mathbb{T}^d$, according to the formula of integral by parts,
\begin{align*}
\sum_{z\sim y}g_t^x(z)-(d+1)g_t^x(y)&=\int_0^{+\infty}e^{-\frac{1}{\sqrt{t}}s}\Big(\sum_{z\sim y}p_s(x,z)-(d+1)p_s(x,y)\Big)ds\\
&=\int_0^{+\infty}e^{-\frac{1}{\sqrt{t}}s}\frac{d}{ds}p_s(x,y)ds\\
&=e^{-\frac{1}{\sqrt{t}}s}p_s(x,y)\Big|_0^{+\infty}+\frac{1}{\sqrt{t}}\int_0^{+\infty}e^{-\frac{1}{\sqrt{t}}s}p_s(x,y)ds\\
&=-1_{\{y=x\}}+\frac{1}{\sqrt{t}}g_t^x(y).
\end{align*}
As a result,
\[
\mathcal{L}G_t^x(\eta)=\sum_{y\in \mathbb{T}^d}(\eta(y)-p)\Big(-1_{\{y=x\}}+\frac{1}{\sqrt{t}}g_t^x(y)\Big)=\frac{1}{\sqrt{t}}G_t^x(\eta)-(\eta(x)-p).
\]

\qed

\subsection{The proof of Equation \eqref{equ 3.2}}\label{appendix 3}

\proof[Proof of Equation \eqref{equ 3.2}]

According to the definition of $g_t^x, g_t^w$,
\begin{align*}
&\sum_{y\in \mathbb{T}^d}\sum_{z\sim y}\left(g_t^x(y)-g_t^x(z)\right)\left(g_t^w(y)-g_t^w(z)\right)\\
&=2(d+1)\sum_{y\in \mathbb{T}^d}g_t^x(y)g_t^w(y)-2\sum_{y\in \mathbb{T}^d}\sum_{z\sim y}g_t^x(y)g_t^w(z)\\
&=2\sum_{y\in \mathbb{T}^d}g_t^x(y)\Big((d+1)g_t^w(y)-\sum_{z\sim y}g_t^w(z)\Big)\\
&=2\sum_{y\in \mathbb{T}^d}g_t^x(y)\int_0^{+\infty}e^{-\frac{1}{\sqrt{t}}s}\Big((d+1)p_s(w,y)-\sum_{z\sim y}p_s(w,z)\Big)ds\\
&=2\sum_{y\in \mathbb{T}^d}g_t^x(y)\int_0^{+\infty}e^{-\frac{1}{\sqrt{t}}s}\Big(-\frac{d}{ds}p_s(w,y)\Big)ds\\
&=2\sum_{y\in \mathbb{T}^d}g_t^x(y)\left(1_{\{y=w\}}-\frac{1}{\sqrt{t}}g_t^x(y)\right)=2g_t^x(w)-\frac{2}{\sqrt{t}}\sum_{y\in \mathbb{T}^d}\left(g_t^x(y)\right)^2\\
&=2g_t^x(w)-\frac{2}{\sqrt{t}}\int_0^{+\infty}\int_0^{+\infty}e^{-\frac{1}{\sqrt{t}}(s+u)}p_{s+u}(x,x)duds.
\end{align*}
By Equation \eqref{equ 2.1 heat estimation},
\[
\int_0^{+\infty}\int_0^{+\infty}p_{s+u}(x,x)duds<+\infty.
\]
Hence,
\begin{align*}
\lim_{t\rightarrow+\infty}\sum_{y\in \mathbb{T}^d}\sum_{z\sim y}\left(g_t^x(y)-g_t^x(z)\right)\left(g_t^w(y)-g_t^w(z)\right)
=2\lim_{t\rightarrow+\infty}g_t^x(w)=2\int_0^{+\infty}p_s(x,w)ds.
\end{align*}

\qed

\subsection{The proof of Lemma \ref{lemma 5.2}}\label{appendix 4}

\proof[Proof of Lemma \ref{lemma 5.2}]

We only need to show that, for any $b\in \mathbb{R}^m$, $b\cdot \frac{1}{a_t}\Lambda_t^{\{x_i\}_{i=1}^m}$ converges in $\hat{\mathbb{P}}_{c,t}^{\{x_i\}_{i=1}^m}$-probability to $b^T\Gamma_{\{x_j\}_{j=1}^m}c$. For any $x\sim y$ and $\eta\in \{0, 1\}^{\mathbb{T}^d}$, we denote by $\mathcal{T}_t^{x,y}(\eta)$ the term
\[
\exp\left\{\frac{a_t}{t}\sum_{i=1}^mc_iG_t^{x_i}(\eta^{x,y})-\frac{a_t}{t}\sum_{i=1}^mc_iG_t^{x_i}(\eta)\right\}.
\]
For $0\leq s\leq t$, we define
\begin{align*}
\mathcal{M}_s^{t, b}=\frac{1}{a_t}\sum_{i=1}^mb_iG_t^{x_i}(\eta_s)&-\frac{1}{a_t}\sum_{i=1}^mb_iG_t^{x_i}(\eta_0)\\
&-\int_0^s\frac{1}{2a_t}\sum_{x\in \mathbb{T}^d}\sum_{y\sim x}\mathcal{T}_t^{x,y}(\eta_u)\left(\sum_{i=1}^mb_i\left(G_t^{x_i}(\eta_u^{x,y})-G_t^{x_i}(\eta_u)\right)\right)du.
\end{align*}
According to Proposition 7.3 in Appendix 1 of \cite{kipnis+landim99}, $\{\mathcal{M}_s^{t, b}\}_{0\leq s\leq t}$ is a martingale. According to the Taylor's expansion formula up to the first order, the term
\[
\int_0^t\frac{1}{2a_t}\sum_{x\in \mathbb{T}^d}\sum_{y\sim x}\mathcal{T}_t^{x,y}(\eta_u)\left(\sum_{i=1}^mb_i\left(G_t^{x_i}(\eta_u^{x,y})-G_t^{x_i}(\eta_u)\right)\right)du
\]
equals ${\rm \uppercase\expandafter{\romannumeral3}}+{\rm \uppercase\expandafter{\romannumeral4}}+o(1)$, where
\[
{\rm \uppercase\expandafter{\romannumeral3}}=\frac{1}{a_t}\int_0^t\sum_{i=1}^mb_i\mathcal{L}G_t^{x_i}(\eta_u)du
\]
and
\begin{align*}
&{\rm \uppercase\expandafter{\romannumeral4}}\\
&=\frac{1}{t}\int_0^t\sum_{x\in \mathbb{T}^d}\sum_{y\sim x}\sum_{i=1}^m\sum_{j=1}^m\left(\eta_u(x)-\eta_u(y)\right)^2c_ib_j
\left(g_t^{x_i}(y)-g_t^{x_i}(x)\right)\left(g_t^{x_j}(y)-g_t^{x_j}(x)\right)du.
\end{align*}
Then, according to Equation \eqref{equ 3.2} and Lemmas \ref{lemma 3.1}, \ref{lemma 3.2},
\[
{\rm \uppercase\expandafter{\romannumeral3}}+{\rm \uppercase\expandafter{\romannumeral4}}=
-b\cdot\frac{1}{a_t}\Lambda_t^{\{x_i\}_{i=1}^m}+b^T\Gamma_{\{x_j\}_{j=1}^m}c+\varepsilon_{4,t},
\]
where $\limsup_{t\rightarrow+\infty}\frac{t}{a_t^2}\log \mathbb{P}_{\nu_p}\left(|\varepsilon_{4,t}|\geq \epsilon\right)=-\infty$ for any $\epsilon>0$.
By Taylor's expansion formula up to the second order and Equation \eqref{equ 3.2}, there exists a constant $M_2<+\infty$ independent of $t$ such that
\begin{equation}\label{equ A.4.1}
\left(\Xi_t^{t,\{x_j\}_{j=1}^m, c}\right)^2\leq \Xi_t^{t,\{x_j\}_{j=1}^m, 2c}\exp\left\{\frac{a_t^2}{t}M_2\right\}
\end{equation}
for sufficiently large $t$. Hence, by Cauchy-Schwarz inequality,
\[
\limsup_{t\rightarrow+\infty}\frac{t}{a_t^2}\log \hat{\mathbb{P}}_{c,t}^{\{x_i\}_{i=1}^m}\left(|\varepsilon_{4,t}|\geq \epsilon\right)=-\infty.
\]
Consequently, since $|G_t^x(\eta)|\leq \sqrt{t}$, we have
\[
\mathcal{M}_t^{t, b}=b\cdot\frac{1}{a_t}\Lambda_t^{\{x_i\}_{i=1}^m}-b^T\Gamma_{\{x_j\}_{j=1}^m}c+\varepsilon_{5,t},
\]
where $\limsup_{t\rightarrow+\infty}\frac{t}{a_t^2}\log \hat{\mathbb{P}}_{c,t}^{\{x_i\}_{i=1}^m}\left(|\varepsilon_{5,t}|\geq \epsilon\right)=-\infty$ for any $\epsilon>0$. Then, according to Doob's inequality, to complete the proof we only need to show that
\[
\lim_{t\rightarrow+\infty}\sum_{0\leq s\leq t}\left(\mathcal{M}_s^{t,b}-\mathcal{M}_{s-}^{t,b}\right)^2=0
\]
in $\hat{\mathbb{P}}_{c,t}^{\{x_i\}_{i=1}^m}$-probability. At each jump moment $s$, $\mathcal{M}_s^{t,b}-\mathcal{M}_{s-}^{t,b}=O(a_t^{-1})$. Then, under $\mathbb{P}_{\nu_p}$, $\sum_{0\leq s\leq t}\left(\mathcal{M}_s^{t,b}-\mathcal{M}_{s-}^{t,b}\right)^2$ is stochastically dominated from above by $\frac{M_3}{a_t^2}\varpi(tM_4)$, where $M_3, M_4<+\infty$ are two constants independent of $t$ and $\{\varpi(t)\}_{t\geq 0}$ is a Poisson process at rate $1$. Therefore, according to Markov inequality, it is easy to check that
\[
\limsup_{t\rightarrow+\infty}\frac{1}{a_t^2}\log\mathbb{P}_{\nu_p}\left(\sum_{0\leq s\leq t}\left(\mathcal{M}_s^{t,b}-\mathcal{M}_{s-}^{t,b}\right)^2\geq \epsilon\right)=-\infty
\]
for any $\epsilon>0$. Then, by Equation \eqref{equ A.4.1} and Cauchy-Schwarz inequality, we have
\[
\limsup_{t\rightarrow+\infty}\frac{1}{a_t^2}\log\hat{\mathbb{P}}_{c,t}^{\{x_i\}_{i=1}^m}\left(\sum_{0\leq s\leq t}\left(\mathcal{M}_s^{t,b}-\mathcal{M}_{s-}^{t,b}\right)^2\geq \epsilon\right)=-\infty
\]
for any $\epsilon>0$ and the proof is complete. 

\qed

\textbf{Acknowledgments.}
The author is grateful to Dr. Linjie Zhao for useful comments.
The author is grateful to financial
supports from the National Natural Science Foundation of China with grant number 12371142 and the Fundamental Research Funds for the Central Universities with grant number 2022JBMC039.

{}
\end{document}